\documentclass[12pt]{article}

\usepackage{amsmath, amssymb}
\usepackage{amsthm}
\usepackage{enumitem}
\usepackage{rotating}

\usepackage{multirow}
\usepackage{float}

\usepackage{hyperref}

\usepackage{color}
\definecolor{darkred}{RGB}{139,0,0}
\definecolor{darkgreen}{RGB}{0,100,0}
\definecolor{darkmagenta}{RGB}{139,0,139}
\definecolor{darkpurple}{RGB}{110,0,180}
\definecolor{darkblue}{RGB}{40,0,200}
\definecolor{darkorange}{RGB}{255,140,0}

\newcommand{\bsx}{\boldsymbol{x}}
\newcommand{\bsh}{\boldsymbol{h}}

\newcommand{\bszero}{\boldsymbol{0}}
\newcommand{\bsgamma}{\boldsymbol{\gamma}}

\newcommand{\bsk}{\boldsymbol{k}}

\newcommand{\bsy}{\boldsymbol{y}}

\newcommand{\rd}{\,{\rm d}}

\newcommand{\RR}{\mathbb{R}}
\newcommand{\NN}{\mathbb{N}}
\newcommand{\ZZ}{\mathbb{Z}}

\newcommand{\cO}{\mathcal{O}}
\newcommand{\cH}{\mathcal{H}}
\newcommand{\cF}{\mathcal{F}}
\newcommand{\EC}{\mathrm{CRI}}
\newcommand{\icomp}{\mathtt{i}}
\newcommand{\abs}[1]{\left\vert#1\right\vert}
\newcommand{\APP}{{\rm APP}}

\newtheorem{definition}{Definition}
\newtheorem{theorem}{Theorem}

\newtheorem{remark}[theorem]{Remark}
\newtheorem{lemma}[theorem]{Lemma}

\begin{document}

\title{Tractability of approximation in the weighted Korobov space in the worst-case setting}

\author{Adrian Ebert, Peter Kritzer, and Friedrich Pillichshammer \thanks{The authors are supported by the Austrian Science Fund (FWF), Projects F5506-N26 (Ebert and Kritzer) and F5509-N26 (Pillichshammer), which are parts of the Special Research Program ``Quasi-Monte Carlo Methods: Theory and Applications'', as well as P34808 (Kritzer).}}

\maketitle

\centerline{\large Dedicated to Pierre L'Ecuyer on the occasion of his $70^{{\rm th}}$ birthday}

\begin{abstract}
In this paper we consider $L_p$-approximation, $p \in \{2,\infty\}$, of periodic functions from weighted Korobov spaces. In particular, we discuss tractability properties of such problems, which means that we aim to relate the dependence of the information complexity on the error demand $\varepsilon$ and the dimension $d$ to the decay rate of the weight sequence $(\gamma_j)_{j \ge 1}$ assigned to the Korobov space. Some results have been well known since the beginning of this millennium, others have been proven quite recently. We give a survey of these findings and will add some new results on the $L_\infty$-approximation problem. To conclude, we give a concise overview of results and collect a number of interesting open problems.
\end{abstract}

\section{Introduction}

In this paper we consider $L_p$-approximation, where $p \in \{2,\infty\}$, of periodic functions from a weighted Korobov space with smoothness parameter $\alpha$ from the viewpoint of Information-Based Complexity. In particular, we study the information complexity $n(\varepsilon,d)$ of these problems, which is the minimal number of information evaluations required to push the approximation error below a certain error demand $\varepsilon \in (0,1)$ for problems in dimension $d \in \NN$. The information classes considered are the class $\Lambda^{{\rm all}}$ consisting of arbitrary continuous linear functionals and the class $\Lambda^{{\rm std}}$ consisting of point evaluations only. Furthermore, we will distinguish between the absolute and the normalized error criterion in the worst-case setting.

If the information complexity $n(\varepsilon,d)$ grows exponentially in $d$ for $d$ tending to infinity, the problem is said to suffer from the curse of dimensionality. Otherwise, for sub-exponential growth rates, the problem is said to be tractable. Initially, only the notions of polynomial and strong polynomial tractability were introduced and studied in the literature. An extensive overview of tractability of multivariate problems can be found in the trilogy \cite{NW08,NW10,NW12}. 

For weighted function classes, one assigns real numbers (weights) to the coordinates in order to model varying influence of 
the single variables on the approximation problem, and one is interested in (matching) necessary and sufficient conditions on the weights which guarantee tractability.  In the particular case of $L_2$-approximation for the weighted Korobov space, matching conditions can be found in the paper \cite{WW99} by Wasilkowski and Wo\'zniakowski for the information class $\Lambda^{{\rm all}}$ and in the paper \cite{NSW04} by Novak, Sloan, and Wo\'zniakowski for $\Lambda^{{\rm std}}$. For $L_\infty$-approximation, results on (strong) polynomial tractability are due to Kuo, Wasilkowski, and Wo\'zniakowski; see \cite{KWW09} for $\Lambda^{{\rm all}}$ and \cite{KWW09b} for $\Lambda^{{\rm std}}$.

After (strong) polynomial tractability, more and finer notions of tractability
of tractability have been introduced with the aim of obtaining a more detailed and clearer picture of the tractability of multivariate problems. Nowadays, there is a variety of finer notions of tractability comprising quasi-polynomial tractability, weak tractability, and uniform weak tractability. The exact definitions will be given in Definition~\ref{def:apf-tract}. Based on this development, many multivariate problems need to be reconsidered in order to classify them further with respect to the newer notions of tractability. This has been done recently in \cite{EP21} for the problem of $L_2$-approximation for weighted Korobov spaces. These results will be summarized in Section~\ref{sec:apf-res_app2}. In the present paper we shall also study the $L_\infty$-case. We derive necessary and sufficient conditions for several notions of tractability (see Section~\ref{sec:apf-appinfty}). The presented conditions are tight, but unfortunately do not match exactly. Here, some problems remain open.

In Section~\ref{sec:apf-conclusion} we give a concise survey of the current state of research in tractability theory of approximation in weighted Korobov spaces and formulate some interesting open questions.

Notation and basic definitions will be introduced in the following section.

\section{Basic definitions} 
\label{sec:apf-basics}

\subsubsection*{Function space setting}

The Korobov space $\cH_{d,\alpha,\bsgamma}$ with weight sequence $\bsgamma=(\gamma_j)_{j \ge 1}$ in $\RR^{+}$ 
is a reproducing kernel Hilbert space with kernel function $K_{d,\alpha,\bsgamma}: [0,1]^d \times [0,1]^d \to \RR$ given by
\begin{equation*}
K_{d,\alpha,\bsgamma}(\bsx,\bsy)
:=
\sum_{\bsh \in \ZZ^d} r_{d,\alpha,\bsgamma}(\bsh) \exp(2 \pi \icomp \bsh \cdot (\bsx-\bsy)),
\end{equation*}
where by ``$\cdot$'' we denote the usual dot product. The corresponding inner product and norm are given by
\begin{equation*}
\langle f,g \rangle_{d,\alpha,\bsgamma}
:=
\sum_{\bsh \in \ZZ^d} \frac1{r_{d,\alpha,\bsgamma}(\bsh)} \, \widehat{f}(\bsh) \, \overline{\widehat{g}(\bsh)}
\quad \text{and} \quad
\|f\|_{d,\alpha,\bsgamma}
=
\sqrt{\langle f,f \rangle_{d,\alpha,\bsgamma}}\
.
\end{equation*}
Here, the Fourier coefficients of a function $f\in\cH_{d,\alpha,\bsgamma}$ are given by
\begin{equation*}
\widehat{f}(\bsh)
=
\int_{[0,1]^d} f(\bsx) \exp(-2\pi \icomp \bsh\cdot \bsx) \rd \bsx,
\end{equation*}
and the decay function equals, for $\bsh=(h_1,\ldots,h_d)$, 
$r_{d,\alpha,\bsgamma}(\bsh) = \prod_{j=1}^d r_{\alpha,\gamma_j}(h_j)$, 
with $\alpha > 1$ (the so-called smoothness parameter of the space), and
\begin{equation*}
r_{\alpha,\gamma}(h)
:=
\left\{\begin{array}{ll}
1 & \text{for } h=0 \, , \\[0.5em]
\gamma /|h|^{\alpha} & \text{for } h \in \ZZ\setminus\{0\}.
\end{array}\right.
\end{equation*}
The kernel $K_{d,\alpha,\bsgamma}$ is well-defined for $\alpha > 1$ and for all $\bsx,\bsy \in [0,1]^d$, since
\begin{equation*}
|K_{d,\alpha,\bsgamma}(\bsx,\bsy)| 
\le 
\sum_{\bsh \in \ZZ^d} r_{d,\alpha,\bsgamma}(\bsh) 
= 
\prod_{j=1}^d \left(1+2 \zeta(\alpha) \gamma_j\right) 
< 
\infty
,
\end{equation*}
where $\zeta$ is the Riemann zeta function (note that $\zeta(\alpha)<\infty$ since $\alpha > 1$).

Furthermore, we assume here that the weights are ordered and satisfy 
\begin{equation*}
	1 \ge \gamma_1 \ge \gamma_2 \ge \cdots > 0.
\end{equation*}

The weighted Korobov space is a popular reference space for quasi-Monte Carlo rules, in particular for lattice rules. See, e.g., \cite[Chapter~4]{LP14} or \cite[Appendix~A]{NW08} and the references therein.

\subsubsection*{Approximation in $\cH_{d,\alpha,\bsgamma}$}

In this paper we consider $L_p$-approximation of functions from the weighted Korobov space $\cH_{d,\alpha,\bsgamma}$ for $p \in \{2,\infty\}$. We consider the operator 
$\APP_{d,p}: \cH_{d,\alpha,\bsgamma} \to L_p([0,1]^d)$ with $\APP_{d,p} (f) = f$ for all $f \in \cH_{d,\alpha,\bsgamma}$. The operator $\APP_{d,p}$ is the embedding from the weighted Korobov space $\cH_{d,\alpha,\bsgamma}$ to the space $L_p([0,1]^d)$.

In order to approximate $\APP_{d,p}$ with respect to the $L_p$-norm $\|\cdot\|_{L_p}$ over $[0,1]^d$, $p \in \{2,\infty\}$, it suffices to employ linear algorithms $A_{n,d}$ that use $n$ information evaluations and are of the form
\begin{equation} \label{eq:apf-alg_form}
A_{n,d}(f)
=
\sum_{i=1}^n T_i(f) \, g_i
\quad \text{for }
f \in \cH_{d,\alpha,\bsgamma}	
\end{equation}
with functions $g_i \in L_p([0,1]^d)$ and bounded linear functionals $T_i \in \cH^\ast_{d,\alpha,\bsgamma}$ for $i = 1,\ldots,n$; see \cite{Ba} 
and also \cite{N16,NW08}. We will assume that the functionals $T_i$ belong to some permissible class of information $\Lambda$. In particular, we study the class $\Lambda^{{\rm all}}$ consisting of the entire dual space $\cH^\ast_{d,\alpha,\bsgamma}$ and the class $\Lambda^{{\rm std}}$, which consists only of point evaluation functionals. Recall that $\cH_{d,\alpha,\bsgamma}$ is a reproducing kernel Hilbert space, which means that point evaluations are continuous linear functionals and therefore $\Lambda^{{\rm std}}$ is a subclass of $\Lambda^{{\rm all}}$. With some abuse of notation we will write $A_{n,d} \in \Lambda$ if $A_{n,d}$ is a linear algorithm of the form \eqref{eq:apf-alg_form} using information from the class $\Lambda$.\\

We remark that in both cases $p=2$ and $p=\infty$, the embedding operator $\APP_{d,p}$ is continuous for all 
$d \in \NN$, which can be seen as follows.
\begin{itemize}
	\item For $p=2$, we have for all $f \in \cH_{d,\alpha,\bsgamma}$ that
	\begin{align*}
	\|\APP_{d,2}(f)\|_{L_2}^2
	&=
	\|f\|_{L_2}^2 
	= \sum_{\bsh \in \ZZ^d} |\widehat{f}(\bsh)|^2
	\\
	&\le 
	\sum_{\bsh \in \ZZ^d} \frac{1}{r_{d,\alpha,\bsgamma}(\bsh)} |\widehat{f}(\bsh)|^2 
	= 
	\|f\|_{d,\alpha,\bsgamma}^2 < \infty
	. 
	\end{align*}
	By considering the choice $f \equiv 1$, it follows that the above inequality is sharp such that the operator norm of $\APP_{d,2}$ is given by
	\begin{equation*}
	\|\APP_{d,2}\| =1.
	\end{equation*}
	\item For $p=\infty$, we have for all $f \in \cH_{d,\alpha,\bsgamma}$ that
	\begin{align*}
		\|\APP_{d,\infty}(f)\|_{L_\infty}
		&=
		\|f\|_{L_\infty}
		=
		\sup_{\bsx \in [0,1]^d} |f(\bsx)|
		=
		\sup_{\bsx \in [0,1]^d} \abs{\langle f,K_{d,\alpha,\bsgamma}(\cdot,\bsx) \rangle_{d,\alpha,\bsgamma}}
		\\
		&\le
		\|f\|_{d,\alpha,\bsgamma} \, \sup_{\bsx \in [0,1]^d}  \|K_{d,\alpha,\bsgamma}(\cdot,\bsx)\|_{d,\alpha,\bsgamma}
		\\
		&=
		\|f\|_{d,\alpha,\bsgamma} \, \sup_{\bsx \in [0,1]^d} \sqrt{K_{d,\alpha,\bsgamma}(\bsx,\bsx)}
		\\
		&=
		\|f\|_{d,\alpha,\bsgamma} \, \left(\sum_{\bsh \in \ZZ^d} r_{d,\alpha,\bsgamma}(\bsh) \right)^{1/2}
		\\
		&=
		\|f\|_{d,\alpha,\bsgamma} \, \left( \prod_{j=1}^d \left(1+2 \zeta(\alpha) \gamma_j\right) \right)^{1/2} 
		< 
		\infty
		.
	\end{align*}
	By considering the choice $f = K_{d,\alpha,\bsgamma}(\cdot,\bsx)$, it follows that the above inequality is sharp such that the operator norm of $\APP_{d,\infty}$ is given by
	\begin{equation*}
		\|\APP_{d,\infty}\| = \left( \prod_{j=1}^d \left(1+2 \zeta(\alpha) \gamma_j\right) \right)^{1/2}
		.
	\end{equation*}
\end{itemize}

\subsubsection*{The worst-case setting}

The worst-case error of an algorithm $A_{n,d}$ as in \eqref{eq:apf-alg_form} is defined as
\begin{equation*}
	e(A_{n,d},\APP_{d,p}) := \sup_{\substack{f \in \cH_{d,\alpha,\bsgamma} \\ \|f\|_{d,\alpha,\bsgamma} \le 1}} \|\APP_{d,p} (f) - A_{n,d}(f)\|_{L_p},
\end{equation*}
and the $n$-th minimal worst-case error with respect to the information class $\Lambda$ is given by
\begin{equation*}
	e(n,\APP_{d,p},\Lambda) := \inf_{A_{n,d} \in \Lambda} e(A_{n,d},\APP_{d,p})
	\, ,
\end{equation*}
where the infimum is extended over all linear algorithms of the form \eqref{eq:apf-alg_form} with information from the class $\Lambda$. In the case $p=\infty$ the essential supremum is used in the calculation of 
$\|\APP_{d,\infty} (f) - A_{n,d}(f)\|_{L_\infty}$.

The initial error, i.e., the error obtained by approximating $f$ by zero, equals 
\begin{align*}
	e(0,\APP_{d,p})	
	&= 
	\sup_{\substack{f \in \cH_{d,\alpha,\bsgamma} \\ \|f\|_{d,\alpha,\bsgamma} \le 1}} \|\APP_{d,p} (f) \|_{L_p}
	\\
	&=
	\|\APP_{d,p}\|	
	= 
	\left\{
	\begin{array}{ll} 
		1 & \mbox{ if $p=2$},\\
		\left( \prod_{j=1}^d \left(1+2 \zeta(\alpha) \gamma_j\right) \right)^{1/2} & \mbox{ if $p=\infty$}.
	\end{array}
	\right.
\end{align*}
Note that for $p=\infty$ the initial error $e(0,\APP_{d,\infty})$ may be exponential in $d$ if it is not properly normalized. In the following analysis, we will therefore consider the normalized as well as the absolute error criterion.

We are interested in how the approximation error of algorithms $A_{n,d}$ depends on the number $n$ of information evaluations used and how it depends on the problem dimension~$d$. To this end, we define the so-called information complexity as
\begin{equation*}
	n(\varepsilon,\APP_{d,p}, \Lambda)
	:=
	\min\{n \in \NN_0 \, : \,  e(n,\APP_{d,p},\Lambda) \le \varepsilon \, \EC_{d,p} \}
\end{equation*}
with $\varepsilon \in (0,1)$ and $d \in \NN$, and where either $\EC_{d,p}=1$ for the absolute error criterion (we then write $n_{\text{abs}}(\varepsilon,\APP_{d,p}, \Lambda)$) and $\EC_{d,p}=e(0,\APP_{d,p})=\|\APP_{d,p}\|$ for the normalized error criterion (then, we write $n_{\text{norm}}(\varepsilon,\APP_{d,p}, \Lambda)$).

\subsubsection*{Useful relations}

In the case of $L_2$-approximation we have $e(0,\APP_{d,2})=1$ and hence the absolute and the normalized error criteria coincide. 
This means that
\begin{equation*}
	n_{\text{norm}}(\varepsilon,\APP_{d,2}, \Lambda) 
	= 
	n_{\text{abs}}(\varepsilon,\APP_{d,2}, \Lambda)
\end{equation*}
and we just write  $n(\varepsilon,\APP_{d,2}, \Lambda)$ for $\Lambda \in \{\Lambda^{{\rm all}}, \Lambda^{{\rm std}} \}$. 

In the case of $L_{\infty}$-approximation the situation is different, since $e(0,\APP_{d,\infty}) > 1$. Hence we only have 
\begin{equation} \label{eq:apf-relnorabs}
	n_{\text{norm}}(\varepsilon,\APP_{d,\infty}, \Lambda) 
	\le 
	n_{\text{abs}}(\varepsilon,\APP_{d,\infty}, \Lambda) \quad \text{for} \quad \Lambda \in \{\Lambda^{{\rm all}}, \Lambda^{{\rm std}} \}
	.
\end{equation}
Furthermore, it is well known, see, e.g., \cite{KPW17}, that $L_2$-approximation is not harder than $L_\infty$-approximation for the absolute error criterion, which means that
for $\Lambda \in \{\Lambda^{{\rm all}}, \Lambda^{{\rm std}} \}$ we have
\begin{equation*}
	n(\varepsilon,\APP_{d,2}, \Lambda) 
	\le 
	n_{\text{abs}}(\varepsilon,\APP_{d,\infty}, \Lambda)
	.
\end{equation*}
Thus, necessary conditions for tractability of $L_2$-approximation in the weighted space $\cH_{d,\alpha,\bsgamma}$ are also necessary conditions for tractability of $L_\infty$-approximation in $\cH_{d,\alpha,\bsgamma}$ for the absolute error criterion.

For the information class $\Lambda^{{\rm all}}$, $L_p$-approximation for $p \in \{2,\infty\}$ can be fully characterized in terms of the eigenvalues of the self-adjoint, compact operator 
\begin{equation*}
	W_d := 	\APP_{d,2}^\ast \APP_{d,2}: \cH_{d,\alpha,\bsgamma} \to \cH_{d,\alpha,\bsgamma}
	.
\end{equation*}

The following well-known lemma (see, e.g., \cite[p.~215]{NW08}) provides information on the eigenpairs of the operator $W_d$. 

\begin{lemma}
	The eigenpairs of the operator $W_d$ are $(r_{d,\alpha,\bsgamma}(\bsk), e_{\bsk})$ with $\bsk \in \ZZ^d$, where for $\bsk \in \ZZ^d$ we set 
	\begin{equation*}
		e_{\bsk}(\bsx) 
		= 
		e_{\bsk,\alpha,\bsgamma} (\bsx)
		:= 
		\sqrt{r_{d,\alpha,\bsgamma}(\bsk)} \, \exp(2 \pi \icomp \bsk\cdot \bsx)\,,\quad \mbox{for $\bsx \in [0,1]^d$.}
	\end{equation*}
\end{lemma}

Furthermore, denote the ordered eigenvalues of $W_d$ by $(\lambda_{d,k})_{k \in \NN}$, where
\begin{equation*}
	\lambda_{d,1} \ge \lambda_{d,2} \ge \lambda_{d,3} \ge \cdots .
\end{equation*}
Note that $\lambda_{d,1}=1$, since $r_{d,\alpha,\bsgamma}(\bszero)=1$ and $\gamma_j \le 1$ for all $j \in \NN$.

We then have the following relations (see, for example, \cite{NW08,TWW} for $p=2$ and \cite[Theorem~2]{KWW09} for $p=\infty$) for the $n$-th minimal error with respect to $\Lambda^{{\rm all}}$, 
\begin{equation*}
e(n,\APP_{d,p},\Lambda^{{\rm all}}) 
=
\left\{\begin{array}{ll}
\lambda_{d,n+1}^{1/2}  & \text{if } p=2, \\[0.5em]
\left( \sum_{k=n+1}^{\infty} \lambda_{d,k} \right)^{1/2} & \text{if } p=\infty.
\end{array}\right.
\end{equation*}
Consequently, 
\begin{equation*}
n(\varepsilon,\APP_{d,2}, \Lambda)=\min \left\{ n : \lambda_{d,n+1} \le \varepsilon^2 \right\}
\end{equation*}
for $p=2$, and
\begin{equation}\label{eq:apf-infcomev}
n(\varepsilon,\APP_{d,\infty}, \Lambda)=\min \left\{ n : \sum_{k=n+1}^{\infty} \lambda_{d,k} \le \varepsilon^2 \,\EC_{d,\infty}^2 \right\}
\end{equation}
for $p=\infty$.

\subsubsection*{Relations to the average-case setting}

Note that \eqref{eq:apf-infcomev} is exactly the same as the information complexity for $L_2$-approximation in the average-case setting for certain spaces (see \cite[p.~190]{NW12} for a general introduction to the average-case setting). Indeed, following the outline in \cite{KSW08}, assume that we are given a sequence of spaces $\cF_d$, $d\in \NN$, and study the operator $\widetilde{\APP}_{d,2}:\cF_d \to L_2([0,1]^d)$ with $\widetilde{\APP}_{d,2} (f)=f$ for $f\in \cF_d$. Furthermore, we assume that $\cF_d$ is equipped with a Gaussian probability measure $\mu_d$, which has mean zero and a covariance function that coincides with the reproducing kernel of the Korobov space $\cH_{d,\alpha,\bsgamma}$, with all parameters as above. I.e., 
\[
 \int_{\cF_d} f(\bsx) f(\bsy) \mu_d (\rd f) = K_{d,\alpha,\bsgamma} (\bsx,\bsy) \quad \forall \bsx,\bsy\in [0,1]^d.
\]
Again, it is of interest to study approximation of $\widetilde{\APP}_{d,2}$ by linear algorithms $A_{n,d}$ of the form \eqref{eq:apf-alg_form}.
The \textit{average-case error} of such an algorithm $A_{n,d}$ is given by 
\[
 e^{\rm avg} (A_{n,d}, \widetilde{\APP}_{d,2}):=\left(\int_{\cF_d} \big\|\widetilde{\APP}_{d,2}(f)-A_{n,d}(f)\big\|_{L_2([0,1]^d)}^2\ \mu_d (\rd f) \right)^2,
\]
and the initial error by
\[
 e^{\rm avg} (0, \widetilde{\APP}_{d,2}):=\left(\int_{\cF_d} \big\|\widetilde{\APP}_{d,2}(f)\big\|_{L_2([0,1]^d)}^2 \ \mu_d (\rd f) \right)^2.
\]
We can also define the $n$-th minimal average-case error of $L_2$-approximation in $\cF_d$ for an information class $\Lambda$ by
\[
 e(n,\widetilde{\APP}_{d,2},\Lambda):=\inf_{A_{n,d}\in\Lambda} e^{\rm avg} (A_{n,d}, \widetilde{\APP}_{d,2}).
\]
Now define, for any Borel set $G$ in $L_2([0,1]^d)$, the inverse image under $\widetilde{\APP}_{d,2}$ by $\widetilde{\APP}_{d,2}^{-1} (G):=\{f\in \cF_d \ : \  \APP_{d,2} (f) \in G\}$ and let $\nu_d :=\mu_d \circ \widetilde{\APP}_{d,2}^{-1}$. Then, $\nu_d$ is a Gaussian measure on $L_2([0,1]^d)$, again with mean zero, 
and a covariance operator $C_{\nu_d}$ given by 
\[
 (C_{\nu_d} f) (\bsx)=\int_{[0,1]^d} K_{d,\alpha,\bsgamma} (\bsx,\bsy) f(\bsy)\rd \bsy \quad \forall \bsx\in [0,1]^d.
\]
For more detailed information we refer to \cite{KSW08} and the references therein.

Using the notation just introduced, there are several relations to be observed between the worst-case setting and the average-case setting. Indeed, it is known that 
the eigenvalues of the covariance operator $C_{\nu_d}$ coincide with the eigenvalues $(\lambda_{d,k})_{k\in\NN}$ of the operator $W_d$ introduced above. 
Furthermore, by making use of the relation between the covariance function of $\mu_d$ and the kernel $K_{d,\alpha,\bsgamma}$, it can easily be shown that 
\[
 e^{\rm avg} (0, \widetilde{\APP}_{d,2}) = \left(\sum_{\bsk\in\ZZ_d} r_{d,\alpha,\bsgamma}(\bsk)\right)^{1/2}=\left(\sum_{k=1}^\infty \lambda_{d,k}\right)^{1/2}.
\]
Hence the initial error of average-case $L_2$-approximation in $\cF_d$ is exactly the same as the initial error of worst-case $L_\infty$-approximation in 
$\cH_{d,\alpha,\bsgamma}$. What is more, if one allows information from $\Lambda^{\rm all}$, we have 
\[
 e(n,\widetilde{\APP}_{d,2},\Lambda^{{\rm all}})=\left(\sum_{k=n+1}^\infty \lambda_{d,k}\right)^{1/2}
\]
for the $n$-th minimal error, i.e., the $n$-th minimal error of average-case $L_2$-approximation in $\cF_d$ equals the $n$-th minimal 
error of worst-case $L_\infty$-approximation in $\cH_{d,\alpha,\bsgamma}$. For the derivation of these results and further details, we refer to \cite[Chapter 6]{TWW}, see also \cite{NW08}.

These observations (which have been pointed out in the literature before) 
imply that the results on $L_\infty$-approximation in $\cH_{d,\alpha,\bsgamma}$ presented here can also be interpreted as
results on average-case $L_2$-approximation in $\cF_d$. Indeed some of the theorems presented on $L_\infty$-approximation below recover some of the results in \cite{KSW08} and the references therein, formulated for the average-case setting there.

\subsubsection*{Notions of tractability} 

An important goal of tractability theory is to analyze which problems suffer from the curse of dimensionality, i.e., whether there exist $C,\delta>0$ such that $n(\varepsilon,\APP_{d,p},\Lambda) \ge C (1+\delta)^d$ for infinitely many $d \in \NN$, and which do not. In the latter case it is then an important task to classify the growth rate of the information complexity with respect to the dimension~$d$ tending to infinity ($d \rightarrow \infty$) and the error threshold~$\varepsilon$ tending to zero ($\varepsilon \rightarrow 0$). Different growth rates are characterized by means of various notions of tractability which are given in the following definition.

\begin{definition}\label{def:apf-tract}
	Consider the approximation problem $\APP_p=(\APP_{d,p})_{d \ge 1}$ for the information class $\Lambda$. We say that for this problem we have:
	\begin{enumerate}[label=\rm{(\alph*)}]
		\item Strong polynomial tractability \textnormal{(SPT)} if there exist non-negative numbers $\tau, C$ such that
		\begin{equation}\label{eq:apf-defSPT}
		n(\varepsilon,\APP_{d,p}, \Lambda) 
		\le 
		C \, \varepsilon^{-\tau} \quad \mbox{for all $d \in \NN$ and all $\varepsilon \in (0,1)$.}
		\end{equation}
		The infimum of all exponents $\tau \ge 0$ such that \eqref{eq:apf-defSPT} holds for some $C \ge 0$ is called the exponent of strong polynomial tractability and is denoted by $\tau^{\ast}(\Lambda)$.
		\item Polynomial tractability \textnormal{(PT)} if there exist non-negative numbers $\tau, \sigma, C$ such that
		\begin{equation*}
		n(\varepsilon,\APP_{d,p}, \Lambda) \le 	C \, \varepsilon^{-\tau} d^\sigma \quad
		\mbox{for all $d \in \NN$ and all $\varepsilon \in (0,1)$.}
		\end{equation*}
		\item Quasi-polynomial tractability \textnormal{(QPT)} if there exist non-negative numbers $t, C$ such that
		\begin{multline}\label{eq:apf-defQPT}
		n(\varepsilon,\APP_{d,p}, \Lambda) 
		\le 
		C \, \exp(t \,(1 + \ln d) (1 + \ln \varepsilon^{-1})) \\ \mbox{for all $d \in \NN$ and all $\varepsilon \in (0,1)$.}
		\end{multline} 
		The infimum of all exponents $t \ge 0$ such that \eqref{eq:apf-defQPT} holds for some $C \ge 0$ is called the exponent of quasi-polynomial tractability and is denoted by $t^{\ast}(\Lambda)$.
		\item Weak tractability \textnormal{(WT)} if
		\begin{equation*}
		\lim_{d + \varepsilon^{-1} \to \infty} \frac{\ln n(\varepsilon,\APP_{d,p}, \Lambda)}{d + \varepsilon^{-1}} = 	0.
		\end{equation*}
		\item $(\sigma,\tau)$-weak tractability \textnormal{($(\sigma,\tau)$-WT)} for positive numbers $\sigma,\tau$ if
		\begin{equation*}
		\lim_{d + \varepsilon^{-1} \to \infty} \frac{\ln n(\varepsilon,\APP_{d,p}, \Lambda)}{d^\sigma + \varepsilon^{-\tau}} = 0.
		\end{equation*}
		\item Uniform weak tractability \textnormal{(UWT)} if $(\sigma,\tau)$-weak tractability holds for all $\sigma,\tau \in (0,1]$.    
	\end{enumerate} 
\end{definition}    

We obviously have the following hierarchy of tractability notions: 
\begin{equation*}
\text{SPT} \Rightarrow \text{PT} \Rightarrow \text{QPT} \Rightarrow \text{UWT} \Rightarrow (\sigma,\tau)\text{-WT}, \quad 
\text{for any choice of } (\sigma,\tau)\in (0,1]^2
.
\end{equation*}
Furthermore, WT coincides with $(\sigma,\tau)$-WT for $(\sigma,\tau)=(1,1)$. \\

The characterization of the applicable tractability classes will be done with respect to decay conditions on the weight sequence $\bsgamma=(\gamma_j)_{j \ge 1}$. To this end, we introduce the following notation.
\begin{itemize}
	\item The infimum of the sequence $\bsgamma$ is denoted by $\bsgamma_I := \inf_{j \ge 1} \gamma_j$.
	\item The sum exponent $s_{\bsgamma}$ is defined as
	\begin{equation*}
	s_{\bsgamma}:=\inf\left\{\kappa>0 \ : \ \sum_{j=1}^{\infty} \gamma_j^{\kappa} < \infty \right\}
	.
	\end{equation*}
	\item The exponent $t_{\bsgamma}$ is defined as 
	\begin{equation*}
	t_{\bsgamma}:=\inf\left\{\kappa>0 \ : \ \limsup_{d \to \infty} \frac{1}{\ln(d+1)} \sum_{j=1}^d \gamma_j^{\kappa} < \infty \right\}
	.
	\end{equation*}
	\item The exponent $u_{\bsgamma,\sigma}$, for $\sigma>0$, is defined as 
	\begin{equation*}
		u_{\bsgamma,\sigma}
		:= 
		\inf\left\{ \kappa >0 \ : \ \lim_{d \rightarrow \infty} \frac{1}{d^{\sigma}} \sum_{j=1}^d \gamma_j^{\kappa} =0\right\}
		.
	\end{equation*} 
\end{itemize}
In the definitions of $s_{\bsgamma}$, $t_{\bsgamma}$, and $u_{\bsgamma,\sigma}$ we use the convention that $\inf \emptyset=\infty$.

\section{The results for $\APP_2$}
\label{sec:apf-res_app2}

A complete overview of necessary and sufficient conditions for tractability of $L_2$-approximation in the weighted Korobov space has recently been published in \cite{EP21}. 

\begin{theorem}\label{thm:L2_all}
	Consider the approximation problem $\APP_2=(\APP_{d,2})_{d \ge 1}$ for the information class 
	$\Lambda^{{\rm all}}$ and let $\alpha>1$. Then we have the following results.
	\begin{enumerate}
		\item (Cf.~\cite{WW99}) Strong polynomial tractability for the class $\Lambda^{\mathrm{all}}$ holds if and only if $s_{\bsgamma}< \infty$. In this case the exponent of strong polynomial tractability is 
		\begin{equation*}
			\tau^{\ast}(\Lambda^{\mathrm{all}}) = 2 \max\left(s_{\bsgamma},\frac1{\alpha}\right).
		\end{equation*}
		\item (Cf.~\cite{WW99}) Strong polynomial tractability and polynomial tractability for the class $\Lambda^{\mathrm{all}}$ are equivalent.
		\item Quasi-polynomial tractability, uniform weak tractability, and weak tractability for the class $\Lambda^{{\rm all}}$ are equivalent and hold if and only if $\bsgamma_I< 1$.
		\item If we have quasi-polynomial tractability, then the exponent of quasi-polynomial tractability satisfies 
		\begin{equation*}
			t^{\ast}(\Lambda^{{\rm all}}) = 2 \max\left(\frac{1}{\alpha} , \frac{1}{\ln \bsgamma_I^{-1}}\right)
			.
		\end{equation*}
		In particular, if $\bsgamma_I=0$, we set $(\ln \bsgamma_I^{-1})^{-1}:=0$ and we have that $t^{\ast}(\Lambda^{{\rm all}}) = \frac{2}{\alpha}$.
		\item For $\sigma >1$, $(\sigma,\tau)$-weak tractability for the class $\Lambda^{\rm{all}}$ holds for all weights $1 \ge \gamma_1 \ge \gamma_2 \ge \cdots > 0$.
	\end{enumerate}
\end{theorem}

\begin{theorem}\label{thm:L2_std}
	Consider multivariate approximation $\APP_2 = (\APP_{d,2})_{d \ge 1}$ for the information class $\Lambda^{{\rm std}}$ and $\alpha>1$. Then we have the following results.
	\begin{enumerate}
		\item (Cf.~\cite{NSW04}) Strong polynomial tractability for the class $\Lambda^{{\rm std}}$ holds if and only if
		\begin{equation*}
			\sum_{j =1}^{\infty} \gamma_j < \infty,
		\end{equation*}
		which implies $s_{\bsgamma} \le 1$. In this case the exponent of strong polynomial tractability satisfies 
		\begin{equation*}
			\tau^{\ast}(\Lambda^{\mathrm{std}}) =2 \max\left(s_{\bsgamma},\frac{1}{\alpha} \right)
			.
		\end{equation*}
		\item (Cf.~\cite{NSW04}) Polynomial tractability for the class $\Lambda^{{\rm std}}$ holds if and only if
		\begin{equation*}	
			\limsup_{d \to \infty} \frac1{\ln (d+1)} \sum_{j=1}^d \gamma_j < \infty.
		\end{equation*}
		\item Polynomial and quasi-polynomial tractability for the class $\Lambda^{{\rm std}}$ are equivalent.
		\item  Weak tractability for the class $\Lambda^{{\rm std}}$ holds if and only if
		\begin{equation*}
		\lim_{d \to \infty} \frac1{d} \sum_{j=1}^d \gamma_j = 0.
		\end{equation*}
		\item For $\sigma \in (0,1]$, $(\sigma,\tau)$-weak tractability for the class $\Lambda^{\rm{std}}$ holds if and only if
		\begin{equation*}
		\lim_{d \to \infty} \frac1{d^\sigma} \sum_{j=1}^d \gamma_j = 0.
		\end{equation*}
		For $\sigma >1$, $(\sigma,\tau)$-weak tractability for the class $\Lambda^{\rm{std}}$ holds for all weights $1 \ge \gamma_1 \ge \gamma_2 \ge \cdots > 0$.
		\item Uniform weak tractability for the class $\Lambda^{\rm{std}}$ holds if and only if
		\begin{equation*}
		\lim_{d \to \infty} \frac1{d^\sigma} \sum_{j=1}^d \gamma_j = 0 \quad \text{for all } \sigma \in (0,1].
		\end{equation*}
	\end{enumerate}
\end{theorem} 

Theorems \ref{thm:L2_all} and \ref{thm:L2_std} imply that in the case of $L_2$-approximation no open questions remain, at least for the currently most common tractability classes.

\section{The results for $\APP_\infty$}
\label{sec:apf-appinfty}

We have the following result for $L_\infty$-approximation in the space $\cH_{d,\alpha,\bsgamma}$.

\begin{theorem}\label{thm:apf-tract_app_infty}
	Consider multivariate approximation $\APP_\infty = (\APP_{d,\infty})_{d \ge 1}$ for the information classes $\Lambda^{{\rm all}}$ and $\Lambda^{{\rm std}}$ for the normalized and absolute error criterion and $\alpha>1$. Then we have the following results.
	\begin{enumerate}
		\item (Cf.~\cite{KWW09} for $\Lambda^{{\rm all}}$ and \cite{KWW09b} for $\Lambda^{{\rm std}}$) The approximation problem is strongly polynomially tractable if and only if $s_{\bsgamma} < 1$. If this holds, then for any $\tau \in (\max(1/\alpha, s_{\bsgamma}), 1)$ we have
		\begin{align*}
			e(n,\APP_{d,\infty},\Lambda^{{\rm all}}) 
			&= 
			\cO\left( n^{-(1-\tau)/(2 \tau)} \right) \quad \mbox{and}
			\\
			e(n,\APP_{d,\infty},\Lambda^{{\rm std}}) 
			&=
			\cO\left( n^{-(1-\tau)/(2 \tau (1 + \tau))} \right)
			,
		\end{align*}
		where in both cases the implied factor is independent of $n$ and $d$.
		\item (Cf.~\cite{KWW09} for $\Lambda^{{\rm all}}$ and \cite{KWW09b} for $\Lambda^{{\rm std}}$) The approximation problem is polynomially tractable if and only if $t_{\bsgamma} < 1$. If this holds, then for any $\tau \in (\max(1/\alpha, t_{\bsgamma}), 1)$ and any $\delta > 0$ we have
		\begin{align*}
			e(n,\APP_{d,\infty},\Lambda^{{\rm all}}) 
			&=
			\cO\left( n^{-(1-\tau)/(2 \tau)} d^{\delta + \zeta(\alpha \tau) t_{\bsgamma} / \tau} \right) \ \ \mbox{ and }
			\\
			e(n,\APP_{d,\infty},\Lambda^{{\rm std}}) 
			&= 
			\cO\left( n^{-(1-\tau)/(2 \tau (1+\tau))} d^{\delta + \zeta(\alpha \tau) t_{\bsgamma} / \tau} \right),
		\end{align*}
		where in both cases the implied factor is independent of $n$ and $d$.
		\item A necessary condition for quasi-polynomial tractability is 
		\begin{equation*}
			\limsup_{d \rightarrow \infty}\frac{1}{\ln (d+1)} \sum_{j=1}^d \gamma_j< \infty,
		\end{equation*}
		which implies $t_{\bsgamma}\le 1$.
		\item A necessary condition for weak tractability is 
		\begin{equation*}
			\lim_{d \rightarrow \infty}\frac{1}{d} \sum_{j=1}^d \gamma_j=0,
		\end{equation*}
		which implies $u_{\bsgamma,1}\le 1$,  and a sufficient condition for weak tractability is $u_{\bsgamma,1}<~\!1$.
		\item A necessary condition for $(\sigma,\tau)$-weak tractability for $\sigma \in (0,1]$ is 
		\begin{equation*}
			\lim_{d \rightarrow \infty}\frac{1}{d^\sigma} \sum_{j=1}^d \gamma_j=0,
		\end{equation*} 
		which implies $u_{\bsgamma,\sigma} \le 1$, and a sufficient condition 
		for $(\sigma,\tau)$-weak tractability is $u_{\bsgamma,\sigma} <1$. \\
		
		For $\sigma >1$, $(\sigma,\tau)$-weak tractability holds for all weights $1 \ge \gamma_1 \ge \gamma_2 \ge \cdots > 0$.
		\item A necessary condition for uniform weak tractability is 
		\begin{equation*}
			\lim_{d \rightarrow \infty}\frac{1}{d^\sigma} \sum_{j=1}^d \gamma_j=0\quad \mbox{for all $\sigma \in (0,1]$},
		\end{equation*}
		which implies $u_{\bsgamma,\sigma} \le 1$ for all $\sigma \in (0,1]$, and a sufficient condition 
		for uniform weak tractability is
		\begin{equation*}
			u_{\bsgamma,\sigma} < 1 \quad \mbox{for all $\sigma \in (0,1]$.}
		\end{equation*} 
	\end{enumerate}
\end{theorem}

\begin{remark}\label{re1}
	Some remarks on Theorem~\ref{thm:apf-tract_app_infty} are in order.
	\begin{enumerate}
		\item So far we only have a necessary condition for QPT, which is 
		\begin{equation}\label{eq:apf-necQPTLinfty2}
			\limsup_{d \rightarrow \infty}\frac{1}{\ln (d+1)} \sum_{j=1}^d \gamma_j< \infty,
		\end{equation}
		and which in turn implies $t_{\bsgamma}\le 1$. However, this condition is very close to the ``if and only if''-condition for PT, which is $t_{\bsgamma}<1$. It is an interesting question whether \eqref{eq:apf-necQPTLinfty2} is already strong enough to imply QPT or whether $t_{\bsgamma}<1$ is really necessary. The latter case would imply that PT and QPT are equivalent.
		\item The necessary and sufficient conditions for the notions of weak tractability in Items~4--6 are very tight, although not 
		matching exactly. How to close these gaps is another interesting problem. Regarding Item 4, we also refer to \cite[Section 6.3]{NW08}, 
		where a corresponding result for $L_2$-approximation in the average-case setting is shown, and this is---as pointed out in our remarks above---equivalent to our result for $L_\infty$-approximation in the worst-case setting. There, the same gap is observed, but the authors of \cite{NW08} point out that at least for general weights the condition $\lim_{d \rightarrow \infty}\frac{1}{d} \sum_{j=1}^d \gamma_j=0$ is not sufficient for weak tractability. Whether a similar observation also holds for the special case of product weights, which are considered in the present paper, remains open.
	\end{enumerate}
\end{remark}

\begin{proof}[of Theorem~\ref{thm:apf-tract_app_infty}]
	Proofs of the results on (strong) polynomial tractability in Items 1 and 2 can be found in \cite[Theorem~11]{KWW09} for the class $\Lambda^{{\rm all}}$ and in \cite[Theorem~11]{KWW09b} for $\Lambda^{{\rm std}}$. \\
	
	Now we consider QPT. From \eqref{eq:apf-infcomev} and the fact that $\lambda_{d,k}\le1$ for all $k\in \NN$, we have for $n=n_{\text{norm}}(\varepsilon,\APP_{d,\infty}, \Lambda^{{\rm all}})$ that
	\begin{equation*}
		\sum_{k=1}^\infty\lambda_{d,k}-n \le \sum_{k=n+1}^\infty\lambda_{d,k} \le \varepsilon^2\,\sum_{k=1}^\infty\lambda_{d,k}
		.
	\end{equation*}
	Hence,
	\begin{equation}\label{eq:apf-complower}
		n\ge(1-\varepsilon^2)\,\sum_{k=1}^\infty\lambda_{d,k}=(1-\varepsilon^2)\,
		\prod_{j=1}^d\left(1+2\zeta(\alpha)\gamma_j\right)
		.
	\end{equation}
	Assume that we have QPT for $L_\infty$-approximation for $\Lambda^{{\rm all}}$ and the normalized error criterion. Then there exist positive $t$ and $C$ such that 
	\begin{equation*}
		C {\rm e}^{t (1+\ln d)(1+\ln \varepsilon^{-1})} \ge n_{\text{norm}}(\varepsilon,\APP_{d,\infty}, \Lambda^{{\rm all}})\ge  (1-\varepsilon^2)\,
		\prod_{j=1}^d\left(1+2\zeta(\alpha)\gamma_j\right)
	\end{equation*}
	for all $d \in \NN$ and all $\varepsilon \in (0,1)$.
	
	Fixing $\varepsilon \in (0,1)$, e.g., choosing $\varepsilon={\rm e}^{-1}$ and taking the logarithm implies the condition 
	\begin{equation*}
		\ln C+  2 t (1+\ln d) 
		\ge 
		\ln \left(\frac{{\rm e}^2-1}{{\rm e}^2}\right) + \sum_{j=1}^d \ln(1+2 \zeta(\alpha) \gamma_j)
	\end{equation*} 
	for all $d\in \NN$. This implies $\lim_{j \rightarrow \infty} \gamma_j=0$. Since $\frac{\ln(1+x)}{x} \rightarrow 1$ for $x \rightarrow 0$, this then implies 
	\begin{equation}\label{eq:apf-cond_Lp_QPTnec}
		\limsup_{d \rightarrow \infty}\frac{1}{\ln (d+1)} \sum_{j=1}^d \gamma_j < \infty .
	\end{equation}
	Thus we have shown that \eqref{eq:apf-cond_Lp_QPTnec} is a necessary condition for QPT for $\Lambda^{{\rm all}}$ and the normalized error criterion. Since QPT for $\Lambda^{{\rm std}}$ implies QPT for $\Lambda^{{\rm all}}$, we find that \eqref{eq:apf-cond_Lp_QPTnec} is also a necessary condition for QPT for $\Lambda^{{\rm std}}$ and the normalized error criterion.
	
	Assume that we have QPT for $L_\infty$-approximation for $\Lambda \in \{\Lambda^{{\rm all}},\Lambda^{{\rm std}}\}$ and the absolute error criterion. Then, according to \eqref{eq:apf-relnorabs}, we have QPT for $L_{\infty}$-approximation for $\Lambda$ and the normalized error criterion, and hence \eqref{eq:apf-cond_Lp_QPTnec} holds. Thus the proof of Item 3 is complete. \\
	
	We now discuss $(\sigma,\tau)$-WT and first consider the necessary conditions. Assume that we have $(\sigma,\tau)$-WT for $\sigma \in (0,1]$ for $L_\infty$-approximation for $\Lambda^{{\rm all}}$ and the normalized error criterion. Then, according to \eqref{eq:apf-complower},
	\begin{align*}
		0 
		&= 
		\lim_{d+\varepsilon^{-1}\rightarrow \infty}\frac{\ln n_{\text{norm}}(\varepsilon,\APP_{d,\infty}, \Lambda^{{\rm all}})}{d^{\sigma}+\varepsilon^{-\tau}} 
		\\
		&\ge
		\lim_{d+\varepsilon^{-1}\rightarrow \infty} \left(\frac{\ln(1-\varepsilon^2)}{d^{\sigma}+\varepsilon^{-\tau}} + \frac{\sum_{j=1}^d \ln(1+2\zeta(\alpha)\gamma_j)}{d^{\sigma}+\varepsilon^{-\tau}}\right)
		.
	\end{align*}
	For fixed $\varepsilon \in (0,1)$ this implies 
	\begin{equation*}
		\lim_{d \rightarrow \infty}\frac{1}{d^\sigma} \sum_{j=1}^d \ln(1+2\zeta(\alpha)\gamma_j)=0
		,
	\end{equation*}
	which in turn implies that 
	\begin{equation}\label{eq:apf-nec_WT}
		\lim_{d \rightarrow \infty}\frac{1}{d^\sigma} \sum_{j=1}^d \gamma_j=0
		.
	\end{equation}
	So \eqref{eq:apf-nec_WT} is a necessary condition for $(\sigma,\tau)$-WT for $\sigma \in (0,1]$ for $\Lambda^{{\rm all}}$ and the normalized error criterion. In the same way as for QPT we see that \eqref{eq:apf-nec_WT} is a necessary condition for $(\sigma,\tau)$-WT for $\sigma \in (0,1]$ for $\Lambda \in \{\Lambda^{{\rm all}},\Lambda^{{\rm std}}\}$ and the normalized and the absolute error criterion. Note that \eqref{eq:apf-nec_WT} implies $u_{\bsgamma,\sigma} \le 1$. This finishes the proof of the necessary conditions in Items 4--6. \\
	
	Next, we discuss sufficient conditions for $(\sigma,\tau)$-WT. In \cite{ZKH} Zeng, Kritzer, and Hickernell constructed a spline algorithm $A_{n,d}^{{\rm spline}}$ based on lattice rules with a prime number $n$ of nodes, for which for arbitrary $\lambda \in (1/2,\alpha/2)$
	\begin{equation}\label{eq:apf-est_ZKH}
		e(A_{n,d}^{{\rm spline}},\APP_{d,\infty}) 
		\le  
		\frac{\sqrt{2}}{n^{\lambda (2\lambda-1)/(4\lambda-1)}} \prod_{j=1}^d 
		\left(1+2^{2\alpha+1} \gamma_j^{1/(2\lambda)}\zeta\left(\frac{\alpha}{2\lambda}\right)\right)^{2\lambda}
		.
	\end{equation}
	Assume that $u_{\bsgamma,\sigma}<1$. Then there exists a $\lambda \in (1/2,\alpha/2)$ such that 
	\begin{equation}\label{eq:apf-bed_stWTabs}
		\lim_{d \rightarrow \infty} \frac{1}{d^{\sigma}}\sum_{j=1}^d \gamma_j^{1/(2\lambda)}=0
		.
	\end{equation}
	We show that \eqref{eq:apf-bed_stWTabs} implies $(\sigma,\tau)$-WT for the class $\Lambda^{{\rm std}}$ and the absolute error criterion (and therefore also for the class $\Lambda^{{\rm all}}$ and, because of \eqref{eq:apf-relnorabs}, the same holds true for the normalized error criterion).
	
	Let 
	\begin{equation*}
		M
		:=
		\left\lceil \left(\frac{\sqrt{2}}{\varepsilon} \prod_{j=1}^d \left(1+2^{2\alpha+1} \gamma_j^{1/(2\lambda)}\zeta\left(\frac{\alpha}{2\lambda}\right)\right)^{2\lambda}\right)^{(4 \lambda-1)/(\lambda(2 \lambda -1))}  \right\rceil
	\end{equation*}
	and let $n$ be the smallest prime number that is greater than or equal to $M$. Note that then, according to Bertrand's postulate,  $n \in [M,2M]$. Hence, according to \eqref{eq:apf-est_ZKH} we have 
	\begin{equation*}
		e(n,\APP_{d,\infty},\Lambda^{{\rm std}}) \le \varepsilon,
	\end{equation*}
	and therefore 
	\begin{multline*}
		n(\varepsilon,\APP_{d,\infty},\Lambda^{{\rm std}}) 
		\le 
		n\\ 
		\le 2 M \le 4   \left(\frac{\sqrt{2}}{\varepsilon} \prod_{j=1}^d 
		\left(1+2^{2\alpha+1} \gamma_j^{1/(2\lambda)}
		\zeta\left(\frac{\alpha}{2\lambda}\right)\right)^{2\lambda}\right)^{(4 \lambda-1)/(\lambda(2 \lambda -1))}
		.
	\end{multline*}
	Taking the logarithm and using that $\ln(1+x) \le x$ for $x\ge 0$ yields
	\begin{equation*}
		\ln n(\varepsilon,\APP_{d,\infty},\Lambda^{{\rm std}}) \le \ln 4 + \frac{4 \lambda-1}{\lambda(2 \lambda -1)} \!\left[\frac{\ln 2}{2} + \ln \varepsilon^{-1} \!+2^{2(\alpha+1)} \lambda \zeta\left(\frac{\alpha}{2\lambda}\right) \sum_{j=1}^d \gamma_j^{1/(2\lambda)}\right]
	\end{equation*}
	and hence 
	\begin{align*}
		\lefteqn{\lim_{d+\varepsilon^{-1}\rightarrow \infty} \frac{\ln n(\varepsilon,\APP_{d,\infty},\Lambda^{{\rm std}})}{d^{\sigma}+\varepsilon^{-\tau}}}\\
		&\le 
		\frac{4 \lambda-1}{\lambda(2 \lambda -1)} \left[
		\lim_{d+\varepsilon^{-1}\rightarrow \infty} \frac{\ln \varepsilon^{-1}}{d^{\sigma}+\varepsilon^{-\tau}} +  2^{2(\alpha+1)}\lambda \zeta\left(\frac{\alpha}{2\lambda}\right) \lim_{d+\varepsilon^{-1}\rightarrow \infty} \frac{1}{d^{\sigma}+\varepsilon^{-\tau}} \sum_{j=1}^d \gamma_j^{1/(2\lambda)}\right] \\
		&= 
		0\,
		,
	\end{align*}
	where we used \eqref{eq:apf-bed_stWTabs}  for the case $\sigma \in (0,1]$ in the last step. If $\sigma >1$ then \eqref{eq:apf-bed_stWTabs} is not required, since $\gamma_j \le 1$ for all $j \in \NN$, and so the limit relation holds anyway. Thus the proof of Items 4--6 is finished.
\end{proof}

\begin{remark}\label{re2}
	Let us briefly comment on the $L_p$-approximation problem $\APP_p=(\APP_{d,p})_{d \ge 1}$ for $p\in (2,\infty)$ and the absolute error criterion. As for the relation between the minimal errors of $L_2$- and $L_{\infty}$-approximation, it can be shown that
	\begin{equation*}
		e(n,\APP_{d,2},\Lambda)
		\le  
		e(n,\APP_{d,p},\Lambda)
		\le  
		e(n,\APP_{d,\infty},\Lambda) 
		\quad \mbox{for all}\ \ n,d\in\NN
		,
	\end{equation*}
	and
	\begin{equation*}
		n(\varepsilon,\APP_{d,2},\Lambda) 
		\le n(\varepsilon,\APP_{d,p},\Lambda) 
		\le n(\varepsilon,\APP_{d,\infty},\Lambda)
		.
	\end{equation*}
	for all $\varepsilon \in (0,1)$, and all $d\in\NN$. \\
	
	Therefore, we can conclude that a sufficient condition on the weights for a certain tractability notion for $L_{\infty}$-approximation is also sufficient for the same tractability notion for the $L_p$-approximation problem. The other way round, every necessary condition on the weights for a certain tractability notion for $L_2$-approximation is also necessary for the same tractability notion for the $L_p$-approximation problem. We summarize the results that are implied by this insight in Table~\ref{tab:apf-tract_overview_Lp} in Section~\ref{sec:apf-conclusion}. However, many of the sufficient and necessary conditions which we obtain in this way are far from matching each other (especially the ones for the class $\Lambda^{{\rm all}}$). Whether a similar observation is also true for the normalized error criterion and for $p\in (2,\infty)$ remains an open question. 
\end{remark}

\section{Overview and formulation of open problems}
\label{sec:apf-conclusion}

In Tables~\ref{tab:apf-tract_overview2}--\ref{tab:apf-tract_overview_Lp} below, we give a concise overview of the known results and conditions for the various tractability notions. Table~\ref{tab:apf-tract_overview2} is concerned with $L_2$-approximation, Table~\ref{tab:apf-tract_overview_infty} with $L_\infty$-approximation, and Table~\ref{tab:apf-tract_overview_Lp} with $L_p$-approximation for $p\in (2,\infty)$.

\subsubsection*{Open problems}

While we have a full picture of the characterizations of the currently most common notions of tractability for the $L_2$-approximation problem, the $L_\infty$-case is only partially solved and several details remain open. In particular, for QPT a sufficient condition is still missing, and for the $(\sigma,\tau)$-weak tractability notions the necessary and sufficient conditions are tight, but do not match (see Remark~\ref{re1}). These cases remain open for the moment.

Furthermore, also the more general $L_p$-approximation problem for arbitrary $p \in (2,\infty)$ remains unsolved. While there are some fragmentary results for the absolute error criterion, the case of the normalized criterion is completely open (see Remark~\ref{re2}).

\begin{table}[H]
	\caption{Overview of the conditions for tractability of the $L_2$-approximation problem $\APP_2$ for product weights satisfying $1 \ge \gamma_1 \ge \gamma_2 \ge \cdots > 0$ (recall that normalized and absolute criterion coincide for $\APP_2$).}
	\label{tab:apf-tract_overview2}
	\centering
	\begin{tabular}{p{4.3cm}p{3.3cm}p{5.8cm}}
		\hline\noalign{\smallskip}
		& $\Lambda^{{\rm all}}$ & $\Lambda^{{\rm std}}$ \\
		\noalign{\smallskip}\hline\noalign{\smallskip}
		{\rm SPT} & $s_{\bsgamma} < \infty$ & $\sum_{j=1}^{\infty} \gamma_j < \infty$ \\[0.65em]
		{\rm PT}   & $s_{\bsgamma} < \infty$ & $\limsup_{d\to\infty} \frac{1}{\ln (d+1)}\sum_{j=1}^d \gamma_j < \infty$ \\[0.65em]
		{\rm QPT} & $ \bsgamma_I < 1$ & $\limsup_{d\to\infty} \frac{1}{\ln (d+1)}\sum_{j=1}^d \gamma_j < \infty$ \\[0.65em]
		{\rm UWT} & $ \bsgamma_I < 1$ & $\lim_{d\to\infty} \frac{1}{d^{\sigma}}\sum_{j=1}^d \gamma_j = 0 \ \forall  \sigma \in (0,1]$ \\[0.65em]
		$(\sigma,\tau)\mbox{-WT}$, $\sigma \in (0,1]$ & $ \bsgamma_I < 1$ &  $\lim_{d\to\infty} \frac{1}{d^{\sigma}}\sum_{j=1}^d 
		\gamma_j = 0$ \\[0.65em]
		\mbox{WT} & $ \bsgamma_I < 1$ & $\lim_{d\to\infty} \frac{1}{d}\sum_{j=1}^d \gamma_j = 0$ \\[0.65em]
		$(\sigma,\tau)\mbox{-WT}$, $\sigma>1$ & no extra condition on $\bsgamma$ &  no extra condition on $\bsgamma$ \\
		\noalign{\smallskip}\hline\noalign{\smallskip}
	\end{tabular}
\end{table}

\begin{table}[H]
	\caption{Overview of the conditions for tractability of the $L_\infty$-approximation problem $\APP_{\infty}$ for product weights satisfying $1 \ge \gamma_1 \ge \gamma_2 \ge \cdots > 0$ (normalized and absolute criterion).}
	\label{tab:apf-tract_overview_infty}
	\centering
	\begin{tabular}{p{4cm}p{7.5cm}}
		\hline\noalign{\smallskip}
		& $\Lambda^{{\rm all}}$ and $\Lambda^{{\rm std}}$ \\
		\noalign{\smallskip}\hline\noalign{\smallskip}
		{\rm SPT} & $s_{\bsgamma} < 1$ \\[0.65em]
		{\rm PT}   & $t_{\bsgamma} < 1$ \\[0.65em]
		{\rm QPT} & nec.:\ \   $\limsup_{d\to\infty}\frac{\sum_{j=1}^d \gamma_j}{\ln (d+1)} < \infty$ \\[0.65em]
		{\rm UWT} & $\left\{\begin{array}{l} \text{nec.:} \ \lim_{d\to\infty}\frac{\sum_{j=1}^d \gamma_j}{d^\sigma} =0 \ \forall  \sigma \in (0,1] \\[0.65em] \text{suff.:} \ u_{\bsgamma,\sigma} <1 \ \forall  \sigma \in (0,1] \end{array}\right.$  \\[1.5em]
		$(\sigma,\tau)\mbox{-WT}$, $\sigma \in (0,1]$ & $\left\{\begin{array}{l} \text{nec.:} \ \lim_{d\to\infty}\frac{\sum_{j=1}^d \gamma_j}{d^\sigma} =0 \\[0.65em] \text{suff.:} \ u_{\bsgamma,\sigma} < 1 \end{array}\right.$ \\[1.5em]
		\mbox{WT} & $\left\{\begin{array}{l} \text{nec.:} \ \lim_{d\to\infty}\frac{\sum_{j=1}^d \gamma_j}{d} =0 \\[0.65em] \text{suff.:} \ u_{\bsgamma,1} < 1 \end{array}\right.$ \\[1.5em]
		$(\sigma,\tau)\mbox{-WT}$, $\sigma>1$ & no extra condition on $\bsgamma$  \\
		\noalign{\smallskip}\hline\noalign{\smallskip}
	\end{tabular}
\end{table}

\begin{table}[H]
	\caption{Overview of the conditions for tractability of the $L_p$-approximation problem $\APP_{p}$, $p \in (2,\infty)$, for product weights satisfying $1 \ge \gamma_1 \ge \gamma_2 \ge \dots > 0$ (absolute criterion).}
	\label{tab:apf-tract_overview_Lp}
	\centering
	\begin{tabular}{p{3.9cm}p{3.6cm}p{5.8cm}}
		\hline\noalign{\smallskip}
		& $\Lambda^{{\rm all}}$ & $\Lambda^{{\rm std}}$ \\
		\noalign{\smallskip}\hline\noalign{\smallskip}
		{\rm SPT}  & $\left\{\begin{array}{l} \text{nec.:} \ s_{\bsgamma}  < \infty \\[0.65em] \text{suff.:} \ s_{\bsgamma} <1 \end{array}\right.$ & $\left\{\begin{array}{l} \text{nec.:} \ s_{\bsgamma}  \le 1 \\[0.65em] \text{suff.:} \ s_{\bsgamma} <1 \end{array}\right.$ \\[1.5em]
		{\rm PT}   & $\left\{\begin{array}{l} \text{nec.:} \ s_{\bsgamma} < 1 \\[0.65em] \text{suff.:} \ t_{\bsgamma}<1 \end{array}\right.$ & $\left\{\begin{array}{l} \text{nec.:} \ \limsup_{d\to\infty}\frac{\sum_{j=1}^d \gamma_j}{\ln (d+1)} < \infty \\[0.65em] \text{suff.:} \ t_{\bsgamma}<1 \end{array}\right.$ \\[1.5em]
		{\rm QPT} & $\left\{\begin{array}{l} \text{nec.:} \ \gamma_I < 1 \\[0.65em] \text{suff.:} \ ? \end{array}\right.$ & $\left\{\begin{array}{l} \text{nec.:} \ \limsup_{d\to\infty}\frac{\sum_{j=1}^d \gamma_j}{\ln (d+1)} < \infty  \\[0.65em] \text{suff.:} \ ? \end{array}\right.$ \\[1.5em]
		{\rm UWT} & $\left\{\begin{array}{l} \text{nec.:} \ \gamma_I < 1   \\[0.65em] \text{suff.:} \ u_{\bsgamma,\sigma} <1 \\[0.65em]\ \ \ \forall  \sigma \in (0,1] \end{array}\right.$ & $\left\{\begin{array}{l} \text{nec.:} \ \lim_{d\to\infty}\frac{\sum_{j=1}^d \gamma_j}{d^\sigma} =0 \\[0.65em]\ \ \ \forall  \sigma \in (0,1] \\[0.65em] \text{suff.:}  \ u_{\bsgamma,\sigma} <1 \ \forall  \sigma \in (0,1] \end{array}\right.$ \\[1.5em]
		$(\sigma,\tau)\mbox{-WT}$, $\sigma \in (0,1]$ & $\left\{\begin{array}{l} \text{nec.:} \ \gamma_I < 1 \\[0.65em] \text{suff.:} \ u_{\bsgamma,\sigma} < 1 \end{array}\right.$ & $\left\{\begin{array}{l} \text{nec.:} \ \lim_{d\to\infty}\frac{\sum_{j=1}^d \gamma_j}{d^\sigma} =0 \\[0.65em] \text{suff.:} \ u_{\bsgamma,\sigma} < 1 \end{array}\right.$ \\[1.5em]
		\mbox{WT} & $\left\{\begin{array}{l} \text{nec.:} \ \gamma_I < 1 \\[0.65em] \text{suff.:} \ u_{\bsgamma,1} < 1 \end{array}\right.$& $\left\{\begin{array}{l} \text{nec.:} \ \lim_{d\to\infty}\frac{\sum_{j=1}^d \gamma_j}{d} =0 \\[0.65em] \text{suff.:} \ u_{\bsgamma,1} < 1 \end{array}\right.$ \\[1.5em]
		$(\sigma,\tau)\mbox{-WT}$, $\sigma>1$ & no extra condition on $\bsgamma$ 	& no extra condition on $\bsgamma$
		\\
		\noalign{\smallskip}\hline\noalign{\smallskip}
	\end{tabular}
\end{table}

\begin{small}
\noindent\textbf{Authors' addresses:}

\medskip
\noindent Adrian Ebert and Peter Kritzer\\
Johann Radon Institute for Computational and Applied Mathematics (RICAM)\\
Austrian Academy of Sciences\\
Altenbergerstr.~69, 4040 Linz, Austria\\
E-mail: \texttt{adrian.ebert@hotmail.com}, \texttt{peter.kritzer@oeaw.ac.at}

\medskip

\noindent Friedrich Pillichshammer\\
Institut f\"{u}r Finanzmathematik und Angewandte Zahlentheorie\\
Johannes Kepler Universit\"{a}t Linz\\
Altenbergerstr.~69, 4040 Linz, Austria\\
E-mail: \texttt{friedrich.pillichshammer@jku.at}
\end{small}

\end{document}